\documentclass{amsart}

\usepackage{graphics}
\usepackage{epsfig}
\usepackage{amsmath}
\usepackage{amssymb,fancyhdr,txfonts,pxfonts}

\usepackage{xcolor}
\usepackage{color, colortbl}
\usepackage{graphicx}
\usepackage{marvosym}

\usepackage{paralist}
\usepackage{multicol}

\usepackage{pstricks,pst-text,pst-grad,pst-node,pst-3dplot,pstricks-add,pst-poly,pst-coil} 
\usepackage{pst-fun,pst-blur} 

\usepackage{hyperref}

\newtheorem{theorem}{Theorem}
\theoremstyle{plain}

\newtheorem{corollary}{Corollary}

\newtheorem{definition}{Definition}

\newtheorem{notation}{Notation}

\newtheorem{proposition}{Proposition}

\begin{document}

\title[Cross Ratio Transforms]{Invariant and Preserving Transforms for Cross Ratio of 4-Points in a line on Desargues Affine Plane}

\author[Orgest ZAKA]{Orgest ZAKA}
\address{Orgest ZAKA: Department of Mathematics-Informatics, Faculty of Economy and Agribusiness, Agricultural University of Tirana, Tirana, Albania}
\email{ozaka@ubt.edu.al, gertizaka@yahoo.com, ozaka@risat.org}

\author[James F. Peters]{James F. Peters}
\address{James F. PETERS: Department of Electrical \& Computer Engineering, University of Manitoba, WPG, MB, R3T 5V6, Canada and Department of Mathematics, Faculty of Arts and Sciences, Ad\.{i}yaman University, 02040 Ad\.{i}yaman, Turkey}
\thanks{The research has been supported by the Natural Sciences \& Engineering Research Council of Canada (NSERC) discovery grant 185986, Instituto Nazionale di Alta Matematica (INdAM) Francesco Severi, Gruppo Nazionale  per le Strutture Algebriche, Geometriche e Loro Applicazioni grant 9 920160 000362, n.prot U 2016/000036 and Scientific and Technological Research Council of Turkey (T\"{U}B\.{I}TAK) Scientific Human Resources Development (BIDEB) under grant no: 2221-1059B211301223.}
\email{James.Peters3@umanitoba.ca}

\dedicatory{Dedicated to Girard Desargues and Karl G. C. von Staudt}


\begin{abstract}
This paper introduces advances in the geometry of the transforms for cross ratio of four points in a line in the Desargues affine plane. The results given here have a clean, based Desargues affine plan axiomatics and definitions of addition and multiplication of points on a line in this plane, and for skew field properties. In this paper are studied, properties and results related to the some transforms for cross ratio for 4-points, in a line, which we divide into two categories, \emph{Invariant} and \emph{Preserving} transforms for cross ratio. The results in this paper are (1) the cross-ratio of four points is \emph{Invariant} under transforms: Inversion, Natural Translation, Natural Dilation, Mobi\"us Transform, in a line of Desargues affine plane. (2) the cross-ratio of four points is \emph{Preserved} under transforms: parallel projection, translations and dilation's in the Desargues affine plane.
\end{abstract}


\maketitle

\section{Introduction and Preliminaries}\label{sec1}

Influenced by the recently achieved results, related to the ratio of 2 and 3 points (see, \cite{ZakaPeters2022DyckFreeGroup}, \cite{ZakaPeters2022InvariantPreserving}), but mainly on the results presented in the paper \cite{CrossRatio2022ZakaPeters} for cross-ratio of four collinear points in a line $\ell^{OI}$, in Desargues affine planes, in this paper we study some transforms regrading to cross-ratio of four collinear points (four point in a line $\ell^{OI}$ on Desargues affine plane). We divide this transforms in two categories \emph{Invariant-Transforms} and \emph{Preserving-Transforms}. 

Earlier, we define addition and multiplication of points in a line on Desarges affine plane, and 
we have prove that on each line on Desargues affine plane, we can construct a skew-field related to these two actions, so $(\ell^{OI}, +, \cdot)$-is a skew-field, this construction has been achieved, simply and constructively, using simple elements of elementary geometry, and only the basic axioms of Desargues affine plane (see \cite{ZakaFilipi2016}, \cite{FilipiZakaJusufi}, \cite{ZakaThesisPhd}, \cite{ZakaPetersIso} ). 
In this paper, we consider dilations and translations entirely in the Desargues affine plane (see \cite{ZakaDilauto}, \cite{ZakaCollineations}, \cite{ZakaThesisPhd}, \cite{ZakaPetersIso}).

The foundations for the study of the connections between axiomatic geometry and algebraic structures were set forth by D. Hilbert \cite{Hilbert1959geometry}. And some classic for this are,  E. Artin \cite{Artin1957GeometricAlgebra}, D.R. Huges and F.C. Piper ~\cite{HugesPiper}, H. S. M Coxeter ~\cite{CoxterIG1969}. 
Marcel Berger in \cite{Berger2009geometry12}, Robin Hartshorne
 in \cite{Hartshorne1967Foundations}, etc. Even earlier, in we works \cite{ZakaDilauto, ZakaFilipi2016, FilipiZakaJusufi, ZakaCollineations, ZakaVertex, ZakaThesisPhd, ZakaPetersIso, ZakaPetersOrder, ZakaMohammedSF, ZakaMohammedEndo} we have brought up quite a few interesting facts about the association of algebraic structures with affine planes and with ’Desargues affine planes’, and vice versa.

In this paper, all results based in geometric intuition, in axiomatic of Desargues affine plane, and in skew-field properties, we utilize a method that is naive and direct, without requiring the concept of coordinates. 


\subsection{Desargues Affine Plane}
Let $\mathcal{P}$ be a nonempty space, $\mathcal{L}$ is a family of
subsets of $\mathcal{P}$. The elements $P$ of $\mathcal{P}$ are points and an element $\ell$ of $\mathcal{L}$ is a line. 

\begin{definition}
The incidence structure $\mathcal{A}=(\mathcal{P}, \mathcal{L},\mathcal{I})$, called affine plane, where satisfies the above axioms:
\begin{description}
\item[1$^o$] For each points $P,Q\in \mathcal{P}$, there is exactly one line $\ell\in \mathcal{L}$ such that $P,Q \in \ell$.

\item[2$^o$] For each point $P\in \mathcal{P}, \ell\in \mathcal{L}, P \not\in \ell$, there is exactly one line $\ell'\in \mathcal{L}$ such that
$P\in \ell'$ and $\ell\cap \ell' = \emptyset$\ (Playfair Parallel Axiom~\cite{Pickert1973PlayfairAxiom}).   Put another way,
if the point $P\not\in \ell$, then there is a unique line $\ell'$ on $P$ missing $\ell$~\cite{Prazmowska2004DemoMathDesparguesAxiom}.

\item[3$^o$] There is a 3-subset of points $\left\{P,Q,R\right\}\subset\mathcal{P}$, which is not a subset of any $\ell$ in the plane.   Put another way,
there exist three non-collinear points $\mathcal{P}$~\cite{Prazmowska2004DemoMathDesparguesAxiom}.
\end{description}
\end{definition}

\begin{figure}[htbp]
	\centering
		\includegraphics[width=0.85\textwidth]{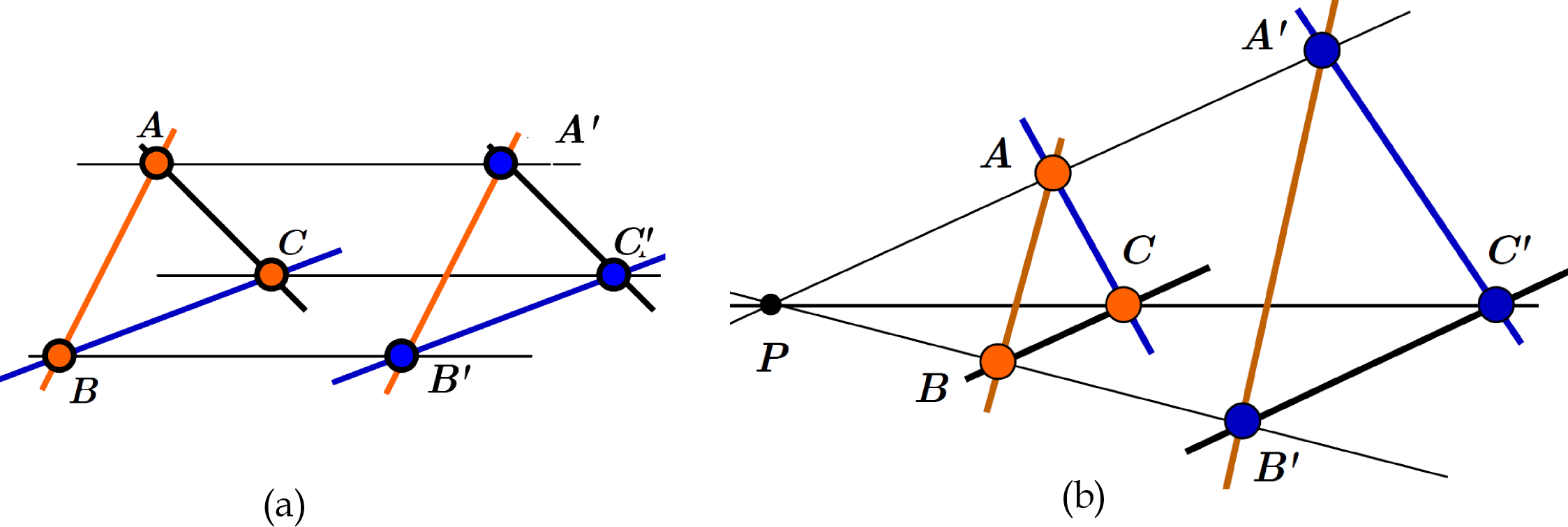}
	\caption{Desargues Axioms: (a) For parallel lines $\ell^{AA_1} \parallel \ell^{BB'} \parallel \ell^{CC'}$; (b) For lines which are cutting in a single point $P$,  $\ell^{AA'} \cap \ell^{BB'} \cap \ell^{CC'}=P$.}
		\label{fig:DesarguesAxiom}
\end{figure}

\emph{\bf Desargues' Axiom, circa 1630}~\cite[\S 3.9, pp. 60-61] {Kryftis2015thesis}~\cite{Szmielew1981DesarguesAxiom}.   Let $A,B,C,A',B',C'\in \mathcal{P}$ and let pairwise distinct lines  $\ell^{AA_1} , \ell^{BB'}, \ell^{CC'}, \ell^{AC}, \ell^{A'C'}\in \mathcal{L}$ such that
\begin{align*}
\ell^{AA'} \parallel \ell^{BB'} \parallel \ell^{CC'} \ \mbox{(Fig.~\ref{fig:DesarguesAxiom}(a))} &\ \mbox{\textbf{or}}\
\ell^{AA'} \cap \ell^{BB'} \cap \ell^{CC'}=P.
 \mbox{(Fig.~\ref{fig:DesarguesAxiom}(b) )}\\
 \mbox{and}\  \ell^{AB}\parallel \ell^{A'B'}\ &\ \mbox{and}\ \ell^{BC}\parallel \ell^{B'C'}.\\
A,B\in \ell^{AB}, A'B'\in \ell^{A'B'}, A,C\in \ell^{AC},  &\ \mbox{and}\ A'C'\in \ell^{A'C'}, B,C\in \ell^{BC},  B'C'\in \ell^{B'C'}.\\
A\neq C, A'\neq C', &\ \mbox{and}\ \ell^{AB}\neq \ell^{A'B'}, \ell^{BC}\neq \ell^{B'C'}.
\end{align*}

Then $\boldsymbol{\ell^{AC}\parallel \ell^{A'C'}}$.   \qquad \textcolor{blue}{$\blacksquare$}

\noindent A {\bf Desargues affine plane} is an affine plane that satisfies Desargues' Axiom. 
\begin{notation}
Three vertexes $ABC$ and $A'B'C'$, which, fulfilling the conditions of the Desargues Axiom, we call \emph{'Desarguesian'}.
\end{notation}

\subsection{Addition and Multiplication of points in a line of Desargues affine plane} $ \quad $\

\begin{figure}[htbp]
\centering%
\includegraphics[width=0.92\textwidth]{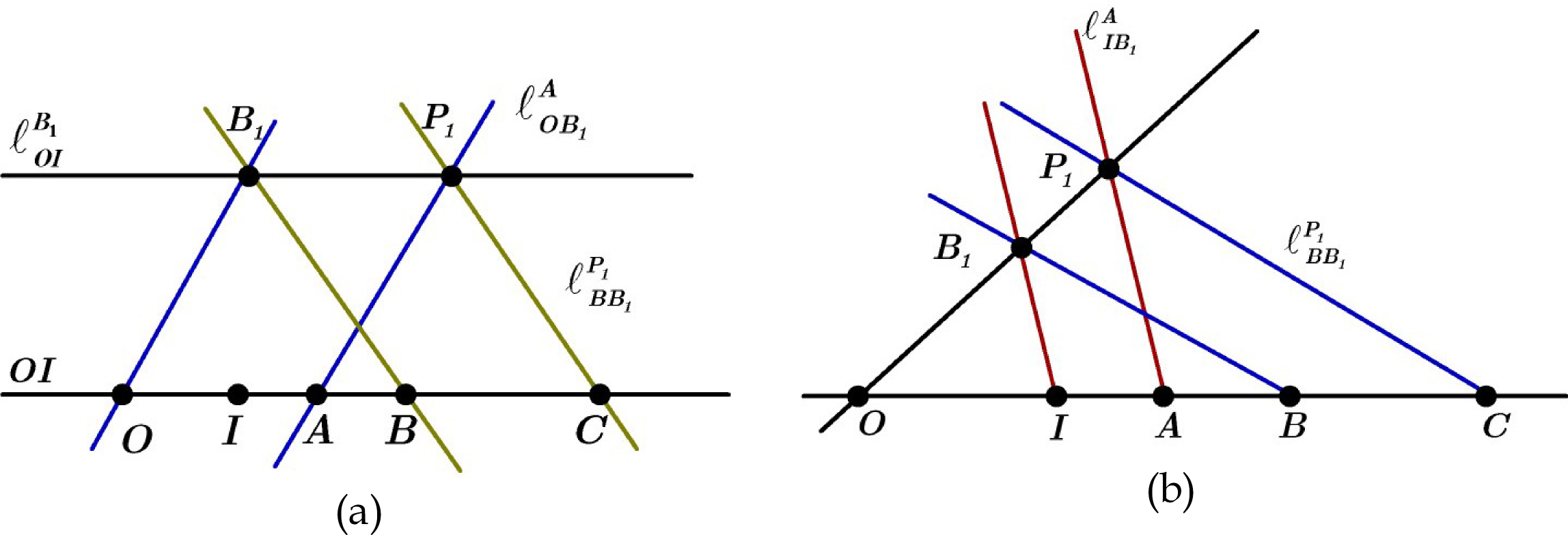}
\caption{ (a) Addition of points in a line in affine plane, 
(b) Multiplication of points in a line in affine plane}
\label{fig:FigureAdMult}
\end{figure}

The process of construct the points $C$ for addition (Figure \ref{fig:FigureAdMult} (a)) and multiplication (Figure \ref{fig:FigureAdMult} (b)) of points in $\ell^{OI}-$line in affine plane, is presented in the tow algorithm form  

\begin{multicols}{2}
\textsc{Addition Algorithm}
\begin{description}
	\item[Step.1] $B_{1}\notin \ell^{OI}$
	\item[Step.2] $\ell_{OI}^{B_{1}}\cap \ell_{OB_{1}}^{A}=P_{1}$
	\item[Step.3] $\ell_{BB_{1}}^{P_{1}}\cap \ell^{OI}=C(=A+B)$
\end{description}

\textsc{Multiplication Algorithm}
\begin{description}
	\item[Step.1] $B_{1}\notin \ell^{OI}$
	\item[Step.2] $\ell_{IB_{1}}^{A}\cap \ell^{OB_{1}}=P_{1}$
	\item[Step.3] $\ell_{BB_{1}}^{P_{1}}\cap \ell^{OI}=C(=A\cdot B)$
\end{description}
\end{multicols}

In \cite{ZakaThesisPhd} and \cite{FilipiZakaJusufi}, we have prove that $(\ell^{OI}, +, \cdot)$ is a skew field in Desargues affine plane, and is field (commutative skew field) in the Papus affine plane.

\begin{definition}
The parallel projection between the two lines in the Desargues affine plane, will be called, a function,
\[P_P : \ell_1 \to  \ell_2, \quad \forall A,B \in \ell_1, \quad AP_P(A) || BP_P(B)
\]
\end{definition}
It is clear that this function is a bijection between any two lines in Desargues affine planes, for this reason, it can also be thought of as isomorphism between two lines.

\begin{definition}
\cite{ZakaCollineations} Dilatation of an affine plane $\mathcal{A}=(\mathcal{P},\mathcal{L}, \mathcal{I})$, called a its collineation $\delta$ such that: $\forall P\neq Q \in \mathcal{P}, \delta{(PQ)}||PQ$.
\end{definition}

\begin{definition}
\cite{ZakaCollineations} Translation of an affine plane $\mathcal{A}=(\mathcal{P},\mathcal{L}, \mathcal{I})$, called identical dilatation $id_{\mathcal{P}}$ his and every other of its dilatation, about which he affine plane has not fixed points.
\end{definition}
\textbf{Some well-known results related to translations and dilation's in Desargues affine planes.}
\begin{itemize}
	\item 
The dilatation set $\textbf{Dil}_{\mathcal{A}}$ of affine plane $\mathcal{A}$ forms a \textbf{group} with respect
to composition $\circ$ (\cite{ZakaCollineations}).
\item The translations set $\textbf{Tr}_{\mathcal{A}}$ of affine plane $\mathcal{A}$ forms a \textbf{group} with respect
to composition $\circ$; which is a sub-group of the dilation group $\left(\textbf{Dil}_{\mathcal{A}}, \circ\right)$ (\cite{ZakaCollineations}).
\item In a affine plane: the group $\left(\textbf{Tr}_{\mathcal{A}}, \circ\right)$ of translations is \textbf{normal
sub-group} of the group of dilatations
$\left(\textbf{Dil}_{\mathcal{A}}, \circ\right)$ (\cite{ZakaCollineations} ).
\item Every dilatation in Desargues affine plane 
	 $\mathcal{A}_{\mathcal{D}}=(\mathcal{P},\mathcal{L}, \mathcal{I})$ which leads a line in it, is an automorphism of skew-fields constructed on the same line $\ell \in \mathcal{L},$ of the plane $\mathcal{A}_{\mathcal{D}}$ (\cite{ZakaDilauto} ).
\item Every translations in Desargues affine plane 
	 $\mathcal{A}_{\mathcal{D}}=(\mathcal{P},\mathcal{L}, \mathcal{I})$ which leads a line in it, is an automorphism of skew-fields constructed on the same line $\ell \in \mathcal{L},$ of the plane $\mathcal{A}_{\mathcal{D}}$ (\cite{ZakaDilauto}).
\item Each dilatation in a Desargues affine plane, $\mathcal{A}_{\mathcal{D}}=(\mathcal{P},\mathcal{L}, \mathcal{I})$ is an isomorphism between skew-fields constructed over isomorphic lines $\ell_1, \ell_2 \in \mathcal{L}$ of that plane (\cite{ZakaPetersIso}).
\item Each translations in a Desargues affine plane, $\mathcal{A}_{\mathcal{D}}=(\mathcal{P},\mathcal{L}, \mathcal{I})$ is an isomorphism between skew-fields constructed over isomorphic lines $\ell_1, \ell_2 \in \mathcal{L}$ of that plane (\cite{ZakaPetersIso}).
\end{itemize}

\subsection{Some algebraic properties of Skew Fields}

I n this section $K$ will denote a skew field~\cite{Herstein1968NR} and $z[K]$ its center, where is the set $K$ such that
$$ z[K]=\left\{k \in K \quad \quad ak=ka, \quad \forall a \in K \right\} .$$

\begin{proposition}
$z[K]$ is a commutative subfield of a skew field $K$.
\end{proposition}

Let now $p\in K$ be a fixed element of the skew field $K$, we will denote by $z_K(p)$ the centralizer in $K$ of the element $p$, where is the set,
\[
z_K(p)=\{k \in K, pk=kp \}.
\]

where $z_K(p)$ is sub skew field of K, but, in general, it is not commutative.

Let $K$ be a skew field, $p\in K$, and let us denote by $[p_K]$ the conjugacy class of $p$:
\[
[p_K]= \left\{q^{-1}pq \quad,\quad q \in K \setminus \{0\} \right\}
\]
If, $p\in z[K]$, for all $q \in K$ we have that $q^{-1}pq=p.$

\subsection{Ratio of two and three points}
In the paper \cite{ZakaPeters2022DyckFreeGroup}, we have done a detailed study, related to the ratio of two and three points in a line of Desargues affine plane. Below we are listing some of the results for ratio of two and three points.
  
\begin{definition} \label{ratio2points}
\cite{ZakaPeters2022DyckFreeGroup} Lets have two different points $A,B \in \ell^{OI}-$line, and $B\neq O$, in Desargues affine plane. We define as ratio of this tow points, a point $R\in \ell^{OI}$, such that,
\[R=B^{-1}A, \qquad \text{
we mark this, with,} \qquad 
R=r(A:B)=B^{-1}A
\]
\end{definition}

For a 'ratio-point' $R \in \ell^{OI}$, and for point $B\neq O$ in line $\ell^{OI}$, is a unique defined point, $A \in \ell^{OI}$, such that $R=B^{-1}A=r(A:B)$.
 
\begin{figure}[htbp]
	\centering
		\includegraphics[width=0.75\textwidth]{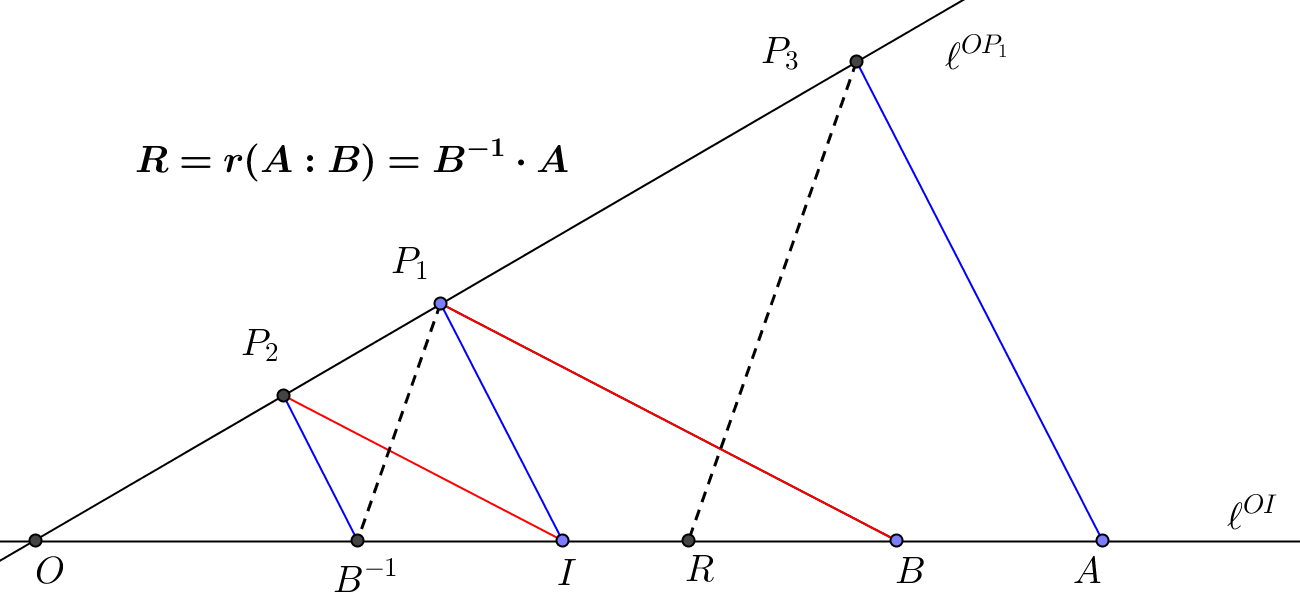}
	\caption{Ilustrate the Ratio-Point, of 2-Points in a line of Desargues affine plane $R=r(A:B)=B^{-1}A$.}
	\label{Ratio2points}
\end{figure}

\textbf{Some results for Ratio of 2-points in Desargues affine plane} (see \cite{ZakaPeters2022DyckFreeGroup}).
\begin{itemize}
\item If have two different points $A,B \in \ell^{OI}-$line, and $B\neq O$, in Desargues affine plane, then, $
r^{-1}(A:B)=r(B:A)$. 
%
\item For three collinear points $A,B,C$ and $C\neq O$, in $\ell^{OI}-$line, have, 
\[
r(A+B:C)=r(A:C)+r(B:C).
\]
\item For three collinear points $A,B,C$ and $C\neq O$, in $\ell^{OI}-$line, have,
\begin{enumerate}
	\item $r(A\cdot B:C)=r(A:C)\cdot B.$
	\item $r(A:B\cdot C)=C^{-1}r(A:C).$
\end{enumerate}
%
\item Let's have the points $A,B$ in the line $\ell^{OI}$ where $B\neq O$.  Then have that, 
\[
r(A:B)=r(B:A) \Leftrightarrow A=B.
\]
%
\item This ratio-map, $r_{B}: \ell^{OI} \to \ell^{OI}$ is a bijection in $\ell^{OI}-$line in Desargues affine plane. 
%
\item The ratio-maps-set $\mathcal{R}_2=\{r_{B}(X),\forall X\in \ell^{OI} \}$, for a fixed point $B$ in the line $\ell^{OI}$, forms a skew-field with 'addition and multiplication' of points. 
This, skew field $(\mathcal{R}_2, +, \cdot)$ is sub-skew field of the skew field $(\ell^{OI}, +, \cdot)$.
\end{itemize}
\textbf{Ratio of three points in a line on Desargues affine plane.} (see \cite{ZakaPeters2022DyckFreeGroup})
\begin{definition}\label{ratiodef}
If $A, B, C$ are three points on a line $\ell^{OI}$ (collinear) in Desargues affine plane, then we define their \textbf{ratio} to be a point $R \in \ell^{OI}$, such that:
\[
(B-C)\cdot R=A-C, \quad \mbox{concisely}\quad R=(B-C)^{-1}(A-C),
\]
and we mark this with  $r(A,B;C)= (B-C)^{-1}(A-C)$.
\end{definition}

\begin{figure}[htbp]
	\centering
		\includegraphics[width=0.9\textwidth]{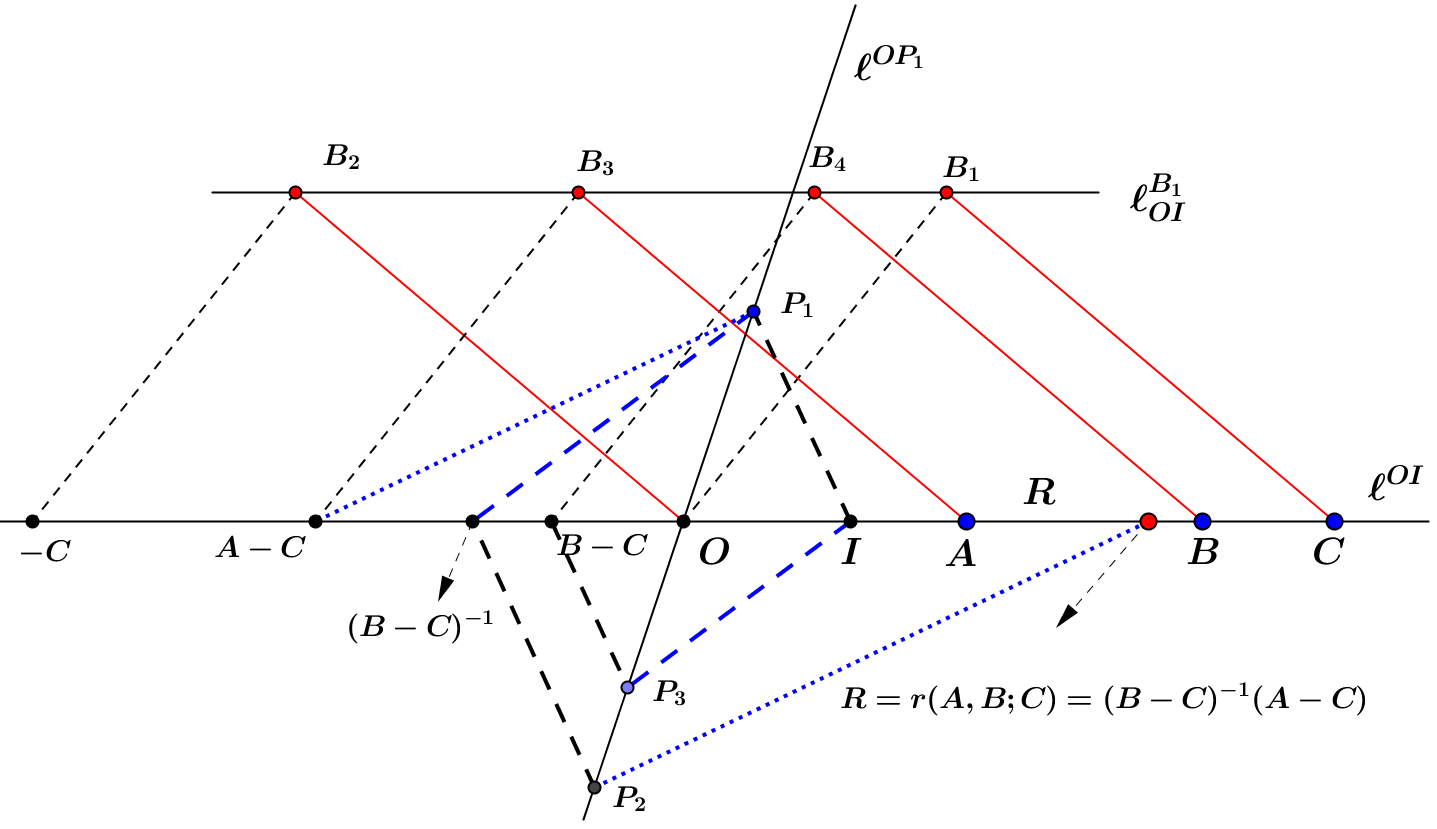}
	\caption{Ratio of 3-Points in a line of Desargues affine plane $R=r(A,B;C)$.}
		\label{ratio3points1}
\end{figure}

\textbf{Some Results for Ratio of 3-points in Desargues affine plane} (\cite{ZakaPeters2022DyckFreeGroup}).
\begin{itemize}
\item \label{reverse.ratio} For 3-points $A,B,C$ in a line $\ell^{OI}$ of Desargues affine plane, we have that,
\[
r(-A,-B;-C)=r(A,B;C).
\]
%
\item \label{inversratio} For 3-points $A,B,C$ in a line $\ell^{OI}$ in the Desargues affine plane, have
\[r^{-1}(A,B;C)=r(B,A;C).
\]
\item If  $A,B,C$, are three different points, and different from point $O$, in a line $\ell^{OI}$ on Desargues affine plane, then
\[r(A^{-1},B^{-1};C^{-1})=B[r(A,B;C)]A^{-1}.\]
%
\item In the Pappus affine plane, for three point different from point $O$, in $\ell^{OI}-$line, we have
$r(A^{-1},B^{-1};C^{-1})=r(A,B;C) \cdot r(B,A;O).$
%
%
\item This ratio-map, $r_{BC}: \ell^{OI} \to \ell^{OI}$ is a bijection in $\ell^{OI}-$line in Desargues affine plane. 
%
%
%
%
\item The ratio-maps-set $\mathcal{R}_3=\{r_{BC}(X),\forall X\in \ell^{OI} \}$, for a different fixed points $B,C$ in $\ell^{OI}-$line, forms a skew-field with 'addition and multiplication' of points in $\ell^{OI}-$line.
This, skew field $(\mathcal{R}_3, +, \cdot)$ is sub-skew field of the skew field $(\ell^{OI}, +, \cdot)$.
\end{itemize}

\subsection{Cross-Ratio in a line of Desargues affine plane}

Let us have the line $\ell^{OI}$ in Desarges affine plane $\mathcal{A_{D}}$, and four points, $A, B, C, D \in \ell^{OI}$

\begin{definition}\label{cross-ratio.def}
If $A, B, C, D$ are four points on a line $\ell^{OI}$ in Desarges affine plane $\mathcal{A_{D}}$, no three of them equal, then we define their cross ratio to be a point:
\[c_r(A,B;C,D)=\left[(A-D)^{-1}(B-D)\right]\left[(B-C)^{-1}(A-C)\right]
\]
\end{definition}

\begin{definition}
If the line $\ell^{OI}$ in Desargues affine plane, is a infinite line (number of points in this line is $+\infty$), we define as follows:
\begin{equation*}
\begin{aligned}
	c_r(\infty, B;C,D) &=(B-D)(B-C)^{-1}\\
	c_r(A,\infty;C,D) &= (A-D)^{-1}(A-C)\\
 c_r(A,B;\infty, D)&=(A-D)^{-1}(B-D) \\
c_r(A,B;C,\infty)&=(B-C)^{-1}(A-C) 
\end{aligned}
\end{equation*}
\end{definition}

From this definition and from ratio definition \ref{ratiodef} we have that,
\begin{itemize}
	\item $c_r(A,B;C,D)=\left[(A-D)^{-1}(B-D)\right]\left[(B-C)^{-1}(A-C)\right]=r(B,A;D) \cdot r(A,B;C).$
	\item $c_r(\infty, B;C,D)=(B-D)(B-C)^{-1}=[(D-B)^{-1}(C-B)]^{-1}=r^{-1}(C,D;B).$
	\item $c_r(A,\infty;C,D)= (A-D)^{-1}(A-C)=(D-A)^{-1}(C-A)=r(C,D;A)$.
	\item $c_r(A,B;\infty, D)=(A-D)^{-1}(B-D)=r(A,B;D)$.
	\item $c_r(A,B;C,\infty)=(B-C)^{-1}(A-C)=r(A,B;C)$.
\end{itemize}

\textbf{Some results for Cross-Ratio of 4-collinear points in Desargues affine plane} (see \cite{CrossRatio2022ZakaPeters}).

\begin{itemize}
	\item If $A,B,C,D$ are distinct points in a $\ell^{OI}-$line, in Desargues affine plane, then
\[ c_r(-A,-B;-C,-D)=c_r(A,B;D,C) \quad\text{and}\quad c_r^{-1}(A,B;C,D)=c_r(A,B;D,C).
\]

\item If $A,B,C,D$ are distinct points in a line, in Desargues affine plane and $I$ is unitary point for multiplications of points in same line, then,
\begin{description}
  \item[(a)] $I-c_r(A,B;C,D)=c_r(A,C;B,D)$
	\item[(b)] $c_r(A,D;B,C)=I-c_r^{-1}(A,B;C,D)$
	\item[(c)] $c_r(A,C;D,B)=[I-c_r(A,B;C,D)]^{-1}$
	\item[(d)] $c_r(A,D;C,B)=[c_r(A,B;C,D)-I]^{-1}c_r(A,B;C,D)$
\end{description}

\item If $A,B,C,D$ are distinct points, and different from zero-point $O$, in a line, in Desargues affine plane and $I$ is unitary point for multiplications of points in same line, have,
\[
c_r(A^{-1},B^{-1}; C^{-1},D^{-1}) = A\cdot c_r(A,B;C,D) \cdot A^{-1}
\]

\item If the point $A\in z[K]$ (center of skew field $K=(\ell^{OI},+,\cdot)$), then, 
\[c_r(A,C;B,D)=c_r(A^{-1},B^{-1}; C^{-1},D^{-1}). \]

\item If $A,B,C,D \in \ell^{OI}$ are distinct points in a line, 
in Desargues affine plane and $I$ is unital point 
for multiplications of points in same line, then equation 
\[
c_r(A,B; C,D) = c_r(B,A;D,C)
\]
it's true, if
\begin{description}
	\item[(a)] points $A,B,C,D$ are in 'center of skew-field' $z[K]$;
	\item[(b)] ratio-points $r(A,B;C)$ are in 'center of skew-field';
	\item[(c)] ratio-point $r(B,A;D)$ are in 'center of skew-field';
	\item[(d)] ratio-point $r(A,B;D)$ is in centaralizer of point $r(A,B;C)$, or vice versa.
\end{description}

\end{itemize}

\section{Some Invariant Transforms for Cross-Ratio of 4-points in a line of Desargues affine plane}

In this section we will see some transformations, for which the cross-ratio are invariant under their action.

\begin{figure}[htbp]
	\centering
		\includegraphics[width=1.0\textwidth]{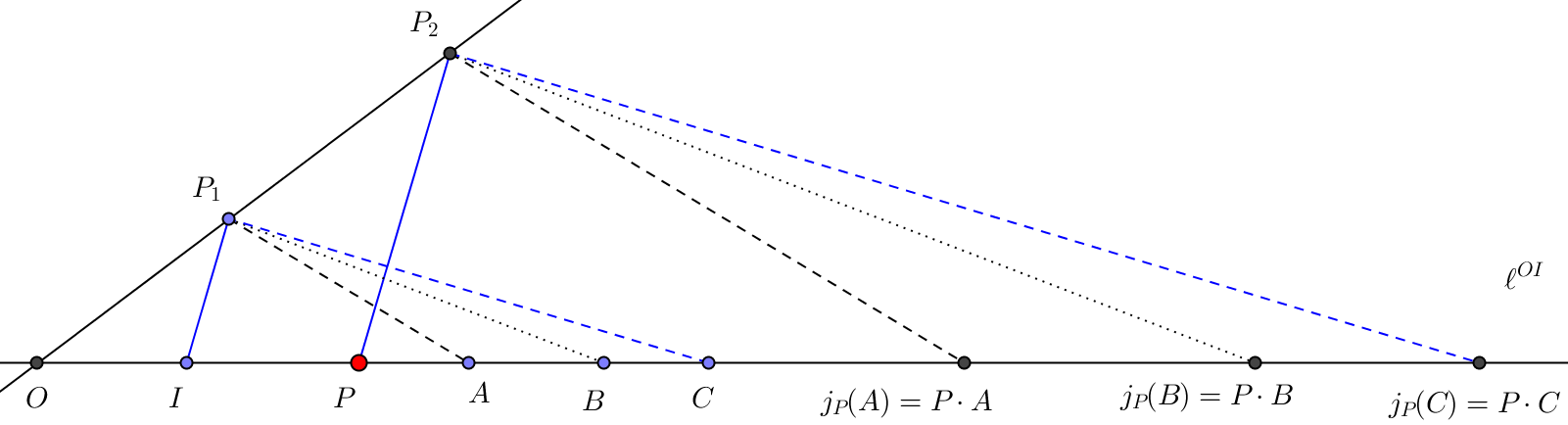}
	\caption{Ilustrate the Inversion of Points, in a line of Desargues affine plane $J_P(A)=P \cdot A$.}
	\label{InversionPoints}
\end{figure}

\begin{figure}[htbp]
	\centering
		\includegraphics[width=1.0\textwidth]{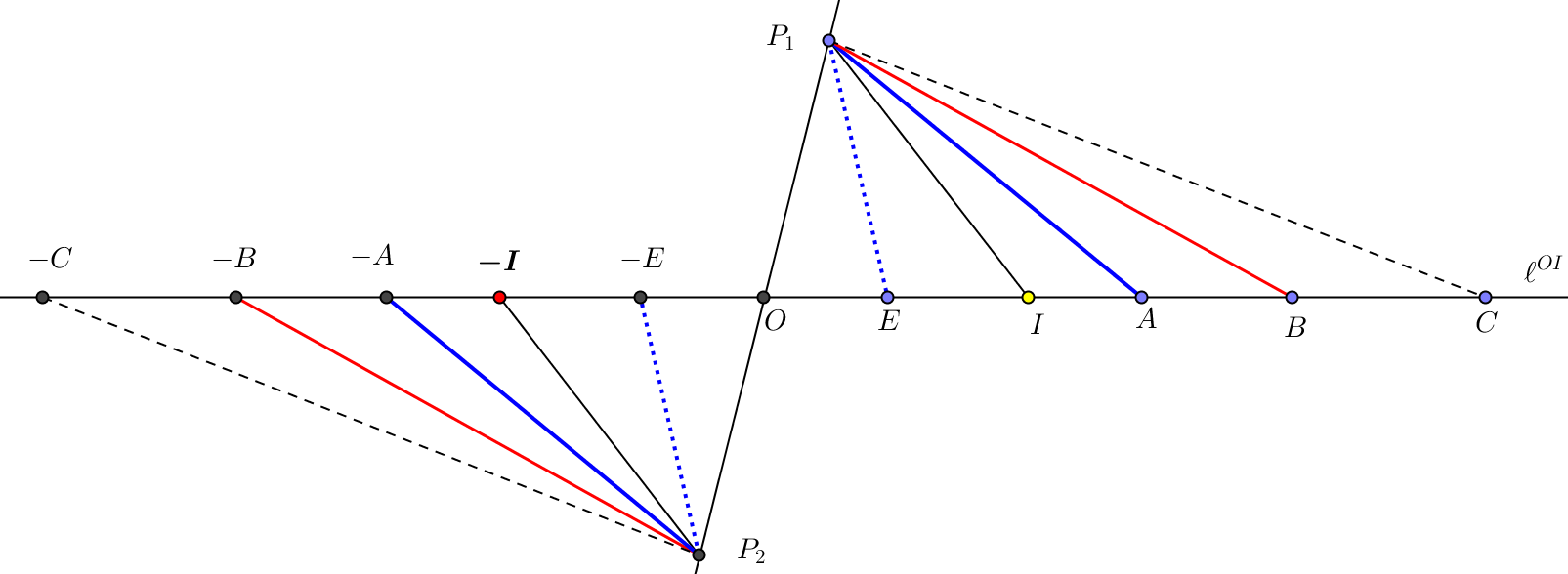}
	\caption{Ilustrate the Reflection of Points, in a line of Desargues affine plane $j_{-I}(A)=-I\cdot A=-A$.}
	\label{ReflectionPoints}
\end{figure}

\begin{figure}[htbp]
	\centering
		\includegraphics[width=1.0\textwidth]{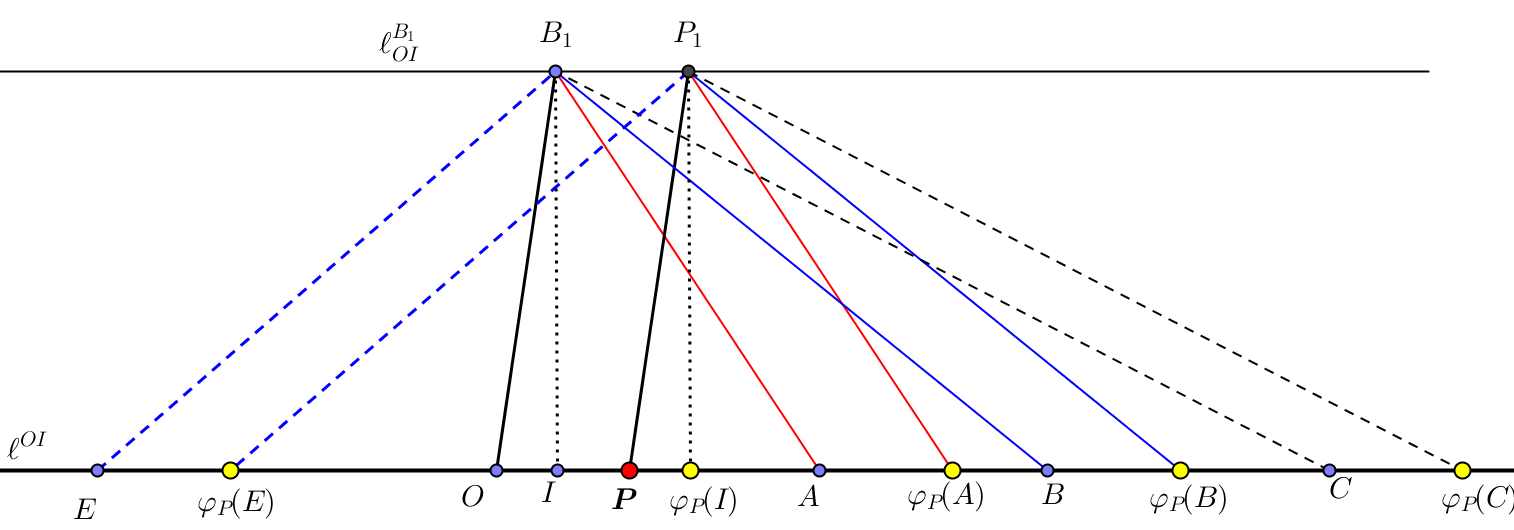}
	\caption{Ilustrate the Natural Translation of Points, in a line of Desargues affine plane $\varphi_P(A)=P+A$.}
	\label{NaturalTrPoints}
\end{figure}

\begin{figure}[htbp]
	\centering
		\includegraphics[width=1.0\textwidth]{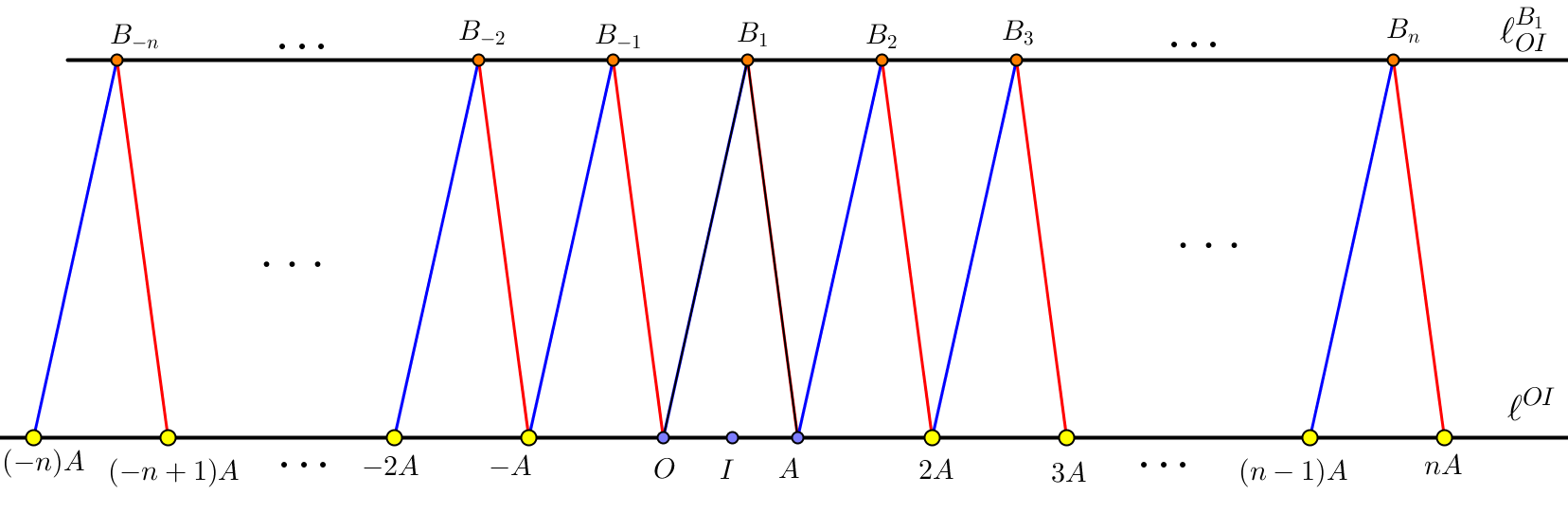}
	\caption{Ilustrate the Natural Translation of Points, in a line of Desargues affine plane $\delta_n(A)=nA$.}
	\label{NaturalDilPoints}
\end{figure}

First, we are defining these transforms and illustrating them with the corresponding figures, respectively: Inversion with Fig.\ref{InversionPoints}, Reflection with Fig.\ref{ReflectionPoints}, Natural Translation with Fig.\ref{NaturalTrPoints} and Natural Dilation with Fig.\ref{NaturalDilPoints}.

\begin{definition}
Inversion of points in $\ell^{OI}-$line, called the map 
\[ j_P:\ell^{OI} \to \ell^{OI}, \]
where $P\in \ell^{OI}$ is fix-point, and which satisfies the condition, 
\[ \forall A \in \ell^{OI} \quad j_P(A)=P \cdot A.\]
\end{definition}

\begin{notation}
Inversion of points in $\ell^{OI}-$line,  
\[ j_P:\ell^{OI} \to \ell^{OI}, \]
where $P=-I \in \ell^{OI}$, called \emph{Involution} or \emph{Reflection about the point $O$ in $\ell^{OI}-$line}
and we have, 
\[ \forall A \in \ell^{OI} \quad j_P(A)=-I \cdot A=-A,\]
where $-A$ is opposite point of point $A$, regrading to addition of points in $\ell^{OI}-$line.
\end{notation}

\begin{definition}
A natural translation with point $P$, of points in $\ell^{OI}-$line, called the map 
\[ \varphi_P:\ell^{OI} \to \ell^{OI}, \]
for a fixed $P \in \ell^{OI}$ which satisfies the condition, 
\[ \forall A \in \ell^{OI} \quad \varphi_P(A)=P+A.\]
\end{definition}

\begin{definition}
A natural Dilation of points in $\ell^{OI}-$line, called the map 
\[ \delta_n:\ell^{OI} \to \ell^{OI}, \]
for a fixed natural number $n \in \mathbb{Z}$ which satisfies the condition, \\
if, $n >0 $, have,
\[ \forall A \in \ell^{OI} \quad \delta_n(A)=nA=\underbrace{A+A+\cdots +A}_{n-times},\]
and if, $n<0$, we have
\[ \forall A \in \ell^{OI} \quad \delta_n(A)=nA=\underbrace{[-A]+[-A]+\cdots +[-A]}_{(-n)-times},\]
where $-A=(-I)\cdot A$ is the oposite point of point $A$, regarding to addition of points in $\ell^{OI}-$line.
\end{definition}

\begin{definition} \label{mobius.transform.def}
Lets have three fixed points $B,C,D\in \ell^{OI}.$ Mobi\"us transform for cross-ratio, we called the map, 
\[
\mu: \ell^{OI} \to \ell^{OI},
\]
which satisfies the condition, 
\[
\forall X \in \ell^{OI}, \quad \mu(X)=c_r(X,B;C,D).
\]
\end{definition}

\begin{theorem}
The cross-ratios are invariant under the natural translation with a point $P$.
\end{theorem}
\proof
From cross-ratio definition \ref{cross-ratio.def}, have,
\[c_r(A,B;C,D)=[(A-D)^{-1}(B-D)][(B-C)^{-1}(A-C)]
\]
so for cross-ratio we have that,

\[
\begin{aligned}
c_r[\varphi_P(A),\varphi_P(B);\varphi_P(C),\varphi_P(D)]
&=c_r(A+P,B+P;C+P,D+P)\\
&= [((A+P)-(D+P))^{-1}((B+P)-(D+P))] \\
& \cdot [((B+P)-(C+P))^{-1}((A+P)-(C+P))]\\
&=[(A+P-D-P)^{-1}(B+P-D-P)] \\
&\cdot [(B+P-C-P)^{-1}(A+P-C-P)]\\
&=[(A-D)^{-1}(B-D)][(B-C)^{-1}(A-C)]\\
&=c_r(A,B;C,D)
\end{aligned}
\]


\begin{theorem}
The Cross-Ratios are invariant under the natural dilation with a fixet $n\in \mathbb{Z}$.
\end{theorem}
\proof
From cross-ratio definition \ref{cross-ratio.def}, we have
\[c_r(A,B;C,D)=[(A-D)^{-1}(B-D)][(B-C)^{-1}(A-C)]
\]
so for cross-ratio of \emph{natural dilation-points} we have that,\\
\textbf{Case.1} For $n>0$, we have,

\[
\begin{aligned}
c_r[\delta_n(A),\delta_n(B);\delta_n(C),\delta_n(D)]
&=c_r(nA,nB;nC,nD)\\
&=[(nA-nD)^{-1}(nB-nD)][(nB-nC)^{-1}(nA-nC)]\\
&=[(n(A-D))^{-1}n(B-D)][(n(B-C))^{-1}n(A-C)]\\
&=[(A-D)^{-1}n^{-1}n(B-D)][(B-C)^{-1}n^{-1}n(A-C)]\\
&=[(A-D)^{-1}(B-D)][(B-C)^{-1}(A-C)]\\
&=c_r(A,B;C,D)
\end{aligned}
\]
\textbf{Case.2} For $n<0$, we mark $m=-n>0$ or $-m=n$ where $m>0$, and have that,
\[
\begin{aligned}
c_r[\delta_n(A),\delta_n(B);\delta_n(C),\delta_n(D)]
&=c_r(nA,nB;nC,nD)\\
&=c_r([-m]A,[-m]B;[-m]C,[-m]D)\\
&=[([-m]A-[-m]D)^{-1}([-m]B-[-m]D)] \\
& \qquad \cdot [([-m]B-[-m]C)^{-1}([-m]A-[-m]C)]\\
&=[(m[-A]-m[-D])^{-1}(m[-B]-m[-D])] \\
& \qquad \qquad \cdot [(m[-B]-m[-C])^{-1}(m[-A]-m[-C])]\\
&=[(m([-A]-[-D]))^{-1}m([-B]-[-D])] \\
& \qquad \qquad \cdot[(m([-B]-[-C]))^{-1}m([-A]-[-C])]\\
&=[([-A]-[-D])^{-1}m^{-1}m([-B]-[-D])] \\
& \qquad \qquad \cdot [([-B]-[-C])^{-1}m^{-1}m([-A]-[-C])]\\
&=[([-A]-[-D])^{-1}([-B]-[-D])] \\
& \qquad \qquad \cdot [([-B]-[-C])^{-1}([-A]-[-C])]\\
&=[([-I][A-D])^{-1}([-I][B-D])] \\
& \qquad \qquad \cdot [([-I][B-C])^{-1}([-I][A-C])]\\
&=[(A-D)^{-1}[-I]^{-1}[-I](B-D)]\\
& \qquad \qquad \cdot [(B-C)^{-1}[-I]^{-1}[-I](A-C)]\\
&=[(A-D)^{-1}([-I][-I])(B-D)] \\
& \qquad \qquad \cdot[(B-C)^{-1}([-I][-I])(A-C)]\\
&=[(A-D)^{-1}(I)(B-D)][(B-C)^{-1}(I)(A-C)]\\
&=[(A-D)^{-1}(B-D)][(B-C)^{-1}(A-C)]\\
&=c_r(A,B;C,D)
\end{aligned}
\]
remember that, $[-I]^{-1}=-I$ and $(-I)\cdot (-I)=I$.

\textbf{Another Proof of Case.2}, for $n<0$, we mark $m=-n>0$ or $-m=n$ where $m>0$, and have that,

\[
c_r[\delta_n(A),\delta_n(B);\delta_n(C),\delta_n(D)]
=c_r(nA,nB;nC,nD)=c_r([-m]A,[-m]B;[-m]C,[-m]D)
\]
so,
\[
c_r([-m]A,[-m]B;[-m]C,[-m]D)=c_r(-[mA],-[mB];-[mC],-[mD])
\]
but from the results listed in section 1, for cross-ratio, we have that,
\[
c_r(-A,-B;-C,-D)=c_r(A,B;C,D), \quad \text{for all different points A,B,C,D} \in \ell^{OI}
\]
therefore
\[
c_r(-[mA],-[mB];-[mC],-[mD])=c_r(mA,mB;mC,mD)
\]
and from Case.1 (since $m>0$), have
\[c_r(mA,mB;mC,mD)=c_r(A,B;C,D).
\]
Hence,
\[ c_r[\delta_n(A),\delta_n(B);\delta_n(C),\delta_n(D)]=c_r(A,B;C,D).
\]
\qed

\begin{theorem}
The Cross-Ratios are invariant under Inversion with a given point $P\in \ell^{OI}$.
\end{theorem}
\proof
From cross-ratio definition \ref{cross-ratio.def}, we have
\[c_r(A,B;C,D)=[(A-D)^{-1}(B-D)][(B-C)^{-1}(A-C)]
\]
so, for cross-ratio of the points $j_P(A),j_P(B);j_P(C),j_P(D) \in \ell^{OI}$  (cross-ratio of Inversion-points) we have that,
\[
\begin{aligned}
c_r[j_P(A),j_P(B);j_P(C),j_P(D)]
&=c_r(PA,PB;PC,PD)\\
&=[(PA-PD)^{-1}(PB-PD)][(PB-PC)^{-1}(PA-PC)]\\
&=[(P(A-D))^{-1}P(B-D)][(P(B-C))^{-1}P(A-C)]\\
&=[(A-D)^{-1}P^{-1}P(B-D)][(B-C)^{-1}P^{-1}P(A-C)] \\
&=[(A-D)^{-1}(P^{-1}P)(B-D)][(B-C)^{-1}(P^{-1}P)(A-C)] \\
&=[(A-D)^{-1}I(B-D)][(B-C)^{-1}I(A-C)]\\
&=[(A-D)^{-1}(B-D)][(B-C)^{-1}(A-C)]\\
&=c_r(A,B;C,D)
\end{aligned}
\]
\qed

\begin{corollary}
The Cross-Ratios are invariant under \emph{reflection about the point $O$ in $\ell^{OI}-$line} in Desargues affine plane.
\end{corollary}
\proof
We have the Inversion with point $(-I)\in \ell^{OI}-$line, and from cross-ratio definition \ref{cross-ratio.def}, we have
\[c_r(A,B;C,D)=[(A-D)^{-1}(B-D)][(B-C)^{-1}(A-C)]
\]
so, for cross-ratio of the points $j_{[-I]}(A),j_{[-I]}(B);j_{[-I]}(C),j_{[-I]}(D) \in \ell^{OI}$  (cross-ratio of Inversion-points) we have that,
\small{
\[
\begin{aligned}
c_r[j_{[-I]}(A),j_{[-I]}(B);j_{[-I]}(C),j_{[-I]}(D)]
&=c_r([-I]A,[-I]B;[-I]C,[-I]D)\\
&=[([-I]A-[-I]D)^{-1}([-I]B-[-I]D)]\\
& \qquad \qquad \cdot[([-I]B-[-I]C)^{-1}([-I]A-[-I]C)]\\
&=[([-I](A-D))^{-1}[-I](B-D)] \\
& \qquad \qquad \cdot[([-I](B-C))^{-1}[-I](A-C)]\\
&=[(A-D)^{-1}[-I]^{-1}[-I](B-D)] \\
& \qquad \qquad \cdot [(B-C)^{-1}[-I]^{-1}[-I](A-C)] \\
&=[(A-D)^{-1}([-I]^{-1}[-I])(B-D)] \\
& \qquad \qquad \cdot [(B-C)^{-1}([-I]^{-1}[-I])(A-C)] \\
&=[(A-D)^{-1}([-I][-I])(B-D)] \\
& \qquad \qquad \cdot [(B-C)^{-1}([-I][-I])(A-C)] \\
&=[(A-D)^{-1}I(B-D)][(B-C)^{-1}I(A-C)]\\
&=[(A-D)^{-1}(B-D)][(B-C)^{-1}(A-C)]\\
&=c_r(A,B;C,D)
\end{aligned}
\]
}
we used the fact that, in the skew field we have true $[-I]^{-1}=-I$ and $[-I][-I]=I$.

Hence
\[
c_r[j_{[-I]}(A),j_{[-I]}(B);j_{[-I]}(C),j_{[-I]}(D)]=c_r(A,B;C,D)
\]


\begin{theorem}
Cross-Ratios are invariant under Mobi\"us transform.
\end{theorem}
\proof
From Mobi\"us transform definition \ref{mobius.transform.def} we have 
\[\mu(X)=c_r(X,B;C,D)=[(X-D)^{-1}(B-D)][(B-C)^{-1}(X-C)]
\]
so, for cross-ratio of the points $\mu(A),\mu(B),\mu(C),\mu(D) \in \ell^{OI}$ first, we calculate, this point, according to following the definition of $\mu-$map, and we have

\begin{itemize}
	\item $\mu(A)=c_r(A,B;C,D)$,
	\item $\mu(B)=c_r(B,B;C,D)$, so
	\[ \begin{aligned}
	\mu(B)&=[(B-D)^{-1}(B-D)][(B-C)^{-1}(B-C)]\\
	&=[I][I]\\
	&=I.
	\end{aligned} \]
Thus
	\[\mu(B)=I \]
	\item $\mu(C)=c_r(C,B;C,D)$, so,
	\[ \begin{aligned}
	\mu(C)&=[(C-D)^{-1}(B-D)][(B-C)^{-1}(C-C)]\\
	&=[(C-D)^{-1}(B-D)][(B-C)^{-1}O]=O\\
	&=[(C-D)^{-1}(B-D)][O]\\
	&=O
	\end{aligned} \]
Thus
	\[\mu(C)=O.\]
	\item $\mu(D)=c_r(D,B;C,D)$, so,
	\[ \begin{aligned}
	\mu(D)&=[(D-D)^{-1}(B-D)][(B-C)^{-1}(D-C)] \\
	&=[O^{-1}(B-D)][(B-C)^{-1}(D-C)]\\
	& \text{(and $O^{-1}=\infty-$(point in infinity))}\\
	&=[\infty][(B-C)^{-1}(D-C)]\\
	&=\infty
\end{aligned} \]
	Thus
	\[ \mu(D)=\infty \]
\end{itemize}

Now, calculate, cross-ratio of points $\mu(A),\mu(B),\mu(C),\mu(D)$, and have
\[
\begin{aligned}
c_r[\mu(A),\mu(B);\mu(C),\mu(D)]
&=c_r(\mu(A),I;O,\infty)\\
&=r(\mu(A),I;O)\\
&=(I-O)^{-1}(\mu(A)-O) \\
&=(I)^{-1}(\mu(A)) \\
&=\mu(A)\\
&=c_r(A,B;C,D)
\end{aligned}
\]
\qed


\section{Transforms which Preserving Cross-Ratios of 4-points in a line of Desargues affine plane}

In this section we prove that the parallel projection, Translations and Dilation of a line in itself or in isomorphic line in Desargues affine plane, \emph{preserving:} the cross-ratios for 4-points. The geometrical interpretations, even in the Euclidean view, are quite beautiful in the above theorems, regardless of the rather rendered figures. This is also the reason that we are giving the proofs in the algebraic sense. So we will always have in mind the close connection of skew field and a line in Desargues affine plane, and the properties of parallel projection, translations and dilation's.

\begin{theorem}
Translations $\varphi$ with trace the same line $\ell^{OI}$ of points $A,B,C,D$, preserve the cross-ratio of this points,
\[
\varphi(c_r(A,B;C,D))=c_r(\varphi(A), \varphi(B);\varphi(C), \varphi(D))
\]
\end{theorem}
\proof
For the proof of that theorem, we refer to the translation 
properties, which are studied in the papers
 (see \cite{ZakaThesisPhd}, \cite{ZakaCollineations}, 
\cite{ZakaDilauto}, \cite{ZakaPetersIso}). 
We will follow algebraic properties for these proofs. 
Whatever the translation $\varphi$, automorphism or isomorphism, 
 $\varphi : \ell^{OI} \rightarrow \ell^{OI}$ or 
$\varphi : \ell^{OI} \rightarrow \ell^{O'I'}$, 
where  $\ell^{OI} \neq \ell^{O'I'}$ and $\ell^{OI} || \ell^{O'I'}$, 
we have true the above equals,

\begin{equation*}
\begin{aligned}
\varphi(c_r(A,B;C,D))&=\varphi\left\{\left[(A-D)^{-1}(B-D)\right]
\left[(B-C)^{-1}(A-C)\right]\right\}\\
&\text{(For this equation, we use the fact that $\varphi$ is homomorphisms)}\\
&=\varphi \left[(A-D)^{-1}(B-D)\right] \cdot \varphi 
\left[(B-C)^{-1}(A-C)\right]\\
&=\varphi[(A-D)^{-1}(B-D)] \cdot \varphi[(B-C)^{-1}(A-C)]\\
&\text{(again, from the fact that $\varphi$ is a homomorphism,)}\\
&=\left\{\varphi[(A-D)^{-1}] \cdot\varphi(B-D) \right\}
\cdot \left\{\varphi[(B-C)^{-1}]\cdot \varphi(A-C)\right\}\\
& \text{(translation $\varphi$ is bijective,)}\\
&=\left\{[\varphi(A-D)]^{-1} \cdot\varphi(B-D) \right\}
\cdot \left\{[\varphi(B-C)]^{-1}\cdot \varphi(A-C)\right\}\\
&\text{(from the fact that $\varphi$ is a homomorphism)}\\
&=\left\{[\varphi(A)-\varphi(D)]^{-1} \cdot [\varphi(B)-\varphi(D) ]\right\}
\cdot \left\{[\varphi(B)-\varphi(C)]^{-1}\cdot [\varphi(A)-\varphi(C)]\right\}\\
&\text{(from cross-ratio definition \ref{cross-ratio.def})}\\
&=c_r(\varphi(A), \varphi(B);\varphi(C), \varphi(D))
\end{aligned}
\end{equation*}
\qed

\begin{theorem}
The parallel projection between the two lines $\ell_1$ and $\ell_2$ in 
Desargues affine plane, \textbf{preserving the cross-ratio} of $4-$points, 
\[
P_P(c_r(A,B;C,D))=c_r(P_P(A), P_P(B);P_P(C),P_P(D))
\]
\end{theorem}
\proof
If $\ell_1 || \ell_2$, we have that the parallel projection is a translation, and have true this theorem.\\
If lines $\ell_1$ and $\ell_2$  they are not parallel (so, they are cut at a single point), we have $A,B,C \in \ell_1$ and $P_P(A), P_P(B), P_P(C), P_P(D) \in \ell_2$.  Also since $P_P$ is a bijection we have that, 

\begin{equation*}
\begin{aligned}
P_P(c_r(A,B;C,D))&=P_P\left\{\left[(A-D)^{-1}(B-D)\right]
\left[(B-C)^{-1}(A-C)\right]\right\}\\
&=P_P[(A-D)^{-1}(B-D)] \cdot P_P[(B-C)^{-1}(A-C)]\\
&=\left\{P_P[(A-D)^{-1}] \cdot P_P(B-D) \right\}
\cdot \left\{P_P[(B-C)^{-1}]\cdot P_P(A-C)\right\}\\
&=\left\{[P_P(A-D)]^{-1} \cdot P_P(B-D) \right\}
\cdot \left\{[P_P(B-C)]^{-1}\cdot P_P(A-C)\right\}\\
&=\left\{[P_P(A)-P_P(D)]^{-1} \cdot [P_P(B)-P_P(D) ]\right\} \\
&\cdot \left\{[P_P(B)-P_P(C)]^{-1}\cdot [P_P(A)-P_P(C)]\right\}\\
&=c_r(\varphi(A), P_P(B);P_P(C), P_P(D))
\end{aligned}
\end{equation*}


\begin{theorem}
Dilation $\delta$ with fixed point in the same line $\ell^{OI}$ of points $A,B,C,D$, preserve the cross-ratio of this points,
\[
\delta(c_r(A,B;C,D))=c_r(\delta(A), \delta(B);\delta(C), \delta(D))
\]
\end{theorem}
\proof
For the proof of that theorem, we refer to the dilation 
properties, which are studied in the papers
 (see \cite{ZakaThesisPhd}, \cite{ZakaCollineations}, 
\cite{ZakaDilauto}, \cite{ZakaPetersIso}). 
We will follow algebraic properties for these proofs. 
Whatever the dilation $\delta$, is a automorphism or a isomorphism, 
 $\delta : \ell^{OI} \rightarrow \ell^{OI}$ or 
$\delta : \ell^{OI} \rightarrow \ell^{O'I'}$, 
where  $\ell^{OI} \neq \ell^{O'I'}$ and $\ell^{OI} || \ell^{O'I'}$, 
we have true the above equals,

\begin{equation*}
\begin{aligned}
\delta(c_r(A,B;C,D))&=\delta\left\{\left[(A-D)^{-1}(B-D)\right]
\left[(B-C)^{-1}(A-C)\right]\right\}\\
&=\delta[(A-D)^{-1}(B-D)] \cdot \delta[(B-C)^{-1}(A-C)]\\
&=\left\{\delta[(A-D)^{-1}] \cdot\delta(B-D) \right\}
\cdot \left\{\delta[(B-C)^{-1}]\cdot \delta(A-C)\right\}\\
&=\left\{[\delta(A-D)]^{-1} \cdot\delta(B-D) \right\}
\cdot \left\{[\delta(B-C)]^{-1}\cdot \delta(A-C)\right\}\\
&=\left\{[\delta(A)-\delta(D)]^{-1} \cdot [\delta(B)-\delta(D) ]\right\}
\cdot \left\{[\delta(B)-\delta(C)]^{-1}\cdot [\delta(A)-\delta(C)]\right\}\\
&=c_r(\delta(A), \delta(B);\delta(C), \delta(D))
\end{aligned}
\end{equation*}
\qed

\section{Declarations}
\subsection*{Funding}
The research by J.F. Peters was supported by the Natural Sciences \& Engineering 
Research Council of Canada (NSERC) discovery grant 185986 and 
Instituto Nazionale di Alta Matematica (INdAM) Francesco Severi, 
Gruppo Nazionale per le Strutture Algebriche, Geometriche e Loro Applicazioni 
grant 9 920160 000362, n.prot U 2016/000036 and Scientific and Technological 
Research Council of Turkey (T\"{U}B\.{I}TAK) Scientific Human Resources 
Development (BIDEB) under grant no: 2221-1059B211301223.

\subsection*{Conflict of Interest Statement}
There is no conflict of interest with any funder.

\subsection*{Authors' contributions}
The authors contribution is as follows: O. Zaka introduced the geometry and its
results in this paper.  J.F. Peters refined and clarified various aspects of
the presented geometry and results in this paper.

\bibliographystyle{amsplain}
\bibliography{RCRrefs}

\end{document}